\documentclass[11pt]{amsart}
\usepackage{amssymb,amscd,amsmath}
\usepackage[all]{xy}
\newtheorem{thm}{Theorem}

\newtheorem{prop}[thm]{Proposition}

\newtheorem{conj}[thm]{Conjecture}

\newtheorem{lem}[thm]{Lemma}
\theoremstyle{remark}

\theoremstyle{definition}

\newcommand{\SO}{\mathrm{SO}}

\newcommand{\R}{{\mathbb{R}}}
\newcommand{\C}{{\mathbb{C}}}
\newcommand{\G}{{\mathbb{G}}}
\renewcommand{\P}{{\mathbb{P}}}

\newcommand{\uB}{{\underline{B}}}
\newcommand{\uC}{{\underline{C}}}
\newcommand{\uG}{{\underline{G}}}

\newcommand{\uN}{{\underline{N}}}

\newcommand{\uT}{{\underline{T}}}

\newcommand{\End}{\mathrm{End}}

\newcommand{\Sym}{\mathrm{Sym}}

\newcommand{\GL}{\mathrm{GL}}
\newcommand{\SL}{\mathrm{SL}}

\newcommand{\Span}{\mathrm{Span}}

\renewcommand{\sl}{\mathfrak{sl}}

\newcommand{\g}{\mathfrak g}

\newcommand{\frakt}{\mathfrak t}

\newcommand{\Alt}{{\raise 2pt\hbox{$\scriptstyle\bigwedge$}}}

\begin{document}
\title{Identities with coefficients in simple compact Lie groups}

\author{Michael Larsen}
\email{mjlarsen@indiana.edu}
\address{Department of Mathematics\\
    Indiana University \\
    Bloomington, IN 47405\\
    U.S.A.}

\author{Aner Shalev}
\email{aner.shalev@mail.huji.ac.il}
\address{Einstein Institute of Mathematics\\
    Hebrew University \\
    Givat Ram, Jerusalem 91904\\
    Israel}


\begin{abstract}
We conjecture that if $G$ is a simple compact Lie group with trivial center, then every $d$-variable non-constant word map with coefficients in $G$ defines a non-constant function on $G^d$.
We prove the conjecture for $A_r$, $B_r$, $E_6$, and $G_2$ using a ping-pong argument.
\end{abstract}

\thanks{ML was partially supported by NSF grant DMS-2001349.
AS was partially supported by ISF grant 700/21 and the
Vinik chair of mathematics which he holds. Both authors acknowledge the support of BSF 
grant 2020/037.}

\maketitle
\begin{center}
\vskip -10pt\noindent
\textit{To Pham Tiep on his 60th birthday.}
\vskip 10pt
\end{center}

Let  $d$ be a positive integer and $w\in F_d$ an element of the free group with $d\ge 1$ generators, which we denote $x_1,\ldots,x_d$.  For any group $G$, evaluation of $w$ at a $d$-tuple of elements of $G$ defines the \emph{word map} $G^d\to G$.  
If $G$ is a Lie group, the word map is smooth; if $\uG$ is an algebraic group over a field $K$, we can regard the word map 
as a morphism of $K$-varieties.  When $\uG$ is semisimple and $w\neq 1$, a theorem of Borel \cite{Borel} asserts that this morphism is always dominant.  In particular, if $\uG$ is defined over $\R$,
and $G := \uG(\R)$, then $G$ is a Lie group, and $w$ is submersive outside a set of measure zero.  From this, it follows that for every semisimple Lie group $G$ and every non-trivial word, the image of the word map has non-empty interior.

By a \emph{word with coefficients in $G$} we mean an element of the free product $G\ast F_d$
for some positive integer $d$.  An element $w\in G\ast F_d$ which is not in
$G$ itself is \emph{non-constant}.
Again, a word with coefficients in $G$ defines an evaluation function $G^d\to G$, which we call the \emph{word map}.  Again, if $G$ is a Lie group, this map is smooth, and if $\uG$ is an algebraic group and $G=\uG(K)$ for some field $K$, then every word with coefficients in $G$ defines a word map which is a morphism of $K$-varieties $\uG^d\to \uG$.

In general, a non-constant word with coefficients in $\uG(K)$ need not determine a dominant word morphism.
For instance, $x_1 g x_1^{-1}$ defines a morphism which factors through $\uG/\uC(g)$, where $\uC(g)$ denotes the centralizer of $g$ in $\uG$.
By the same construction, a non-constant word with coefficients in a Lie group $G$ may define a smooth map $G^d\to G$ whose image has no interior points.
On the other hand, there do exist Lie groups for which all non-constant words define \emph{non-constant} maps.  This paper is motivated by the following conjecture:

\begin{conj}
Let $G$ be a (connected) simple compact Lie group with trivial center.  Then every non-constant word defines a non-constant word map.
\end{conj}

It is easy to see that if $G$ has a non-trivial center, the claim does not hold.  Indeed, let $g$ denote any non-trivial central element of $G$.  Then $x_1 gx_1^{-1}$ is a non-constant word which takes the constant value $g$.
We note also that the conjecture cannot be expected to hold in general for semisimple compact groups with trivial center.  Indeed, if $G = G_1\times G_2$ and $g_1$ and $g_2$ are nontrivial elements of $G_1\times \{1\}$ and $\{1\}\times G_2$
respectively, then $[[x_1,g_1],[x_2,g_2]]$ takes constant value $1$.

We say a word is \emph{normalized} if it evaluates to the identity at $(x_1,\ldots,x_d) = (1,1,\ldots,1)$.  Since $w(1,1, \ldots , 1)^{-1}w(x_1, \ldots , x_d)$ is normalized, it suffices to prove the conjecture for normalized words.
We define a \emph{basic word} in $G\ast F_d$ to be any expression of the form $g x_i g^{-1}$
or $g x_i^{-1} g^{-1}$, where $g\in G$ and $i\in [1,d]$.   We say such a word is of \emph{type} $x_i$ or $x_i^{-1}$ respectively, and we say it has \emph{coefficient} $g$.

\begin{lem}
\label{decomp}
Every normalized word $w\in G\ast F_d$ can be expressed as a product
$w = w_1\cdots w_n$,
where the $w_i$ are basic words.
\end{lem}

\begin{proof}
Every $w$ can be written as a product of elements in the set 
$$G\coprod \{x_1^\pm,\ldots,x_d^\pm\},$$
so we may write 
$$w = g_0 y_1 g_1\cdots y_n g_n,$$
with $g_i\in G$ and $y_i\in  \{x_1^\pm,\ldots,x_d^\pm\}$, and since $w$ is normalized, $g_0g_1\cdots g_n=1$.
For $n=0$, the statement is trivial.  For $n>0$, we write
$$w = (g_0 y_1 g_0^{-1}) w',$$
where either $g_0 y_1 g_0^{-1}$ or $(g_0 y_1 g_0^{-1})^{-1}$ is of the form $g_0 x_j x_0^{-1}$.
Then $w'$ is again normalized, and the lemma holds by induction on $n$.
\end{proof}
Note that we are justified in assuming that no two consecutive basic words whose types are mutually inverse have the same coefficient, since otherwise there is a shorter representation of $w$
as a product of basic words.

If $G$ is a compact connected Lie group, by a well-known theorem of Chevalley,
there exists a connected semisimple algebraic group $\uG$ over $\R$ such that $\uG(\R)$ is isomorphic to $G$ as a topological group.  

\begin{conj}
\label{indep}
Let $\uG$ be a simple algebraic group over $\R$ with $G=\uG(\R)$ compact,  let $\gamma_1,\ldots,\gamma_r\in G$ be elements which are pairwise distinct modulo the center of $G$, and let $d$ be a positive integer.  Then there exist elements $g_1,\ldots,g_d\in G$ such that $\{\gamma_i g_j \gamma_i^{-1}\mid i\in [1,r], j\in [1,d]\}$  generates a free subgroup of $G$ of rank $dr$.
\end{conj}

We remark that there exist elements $\delta,\ldots,\delta_d\in G$ such that all products $\gamma_i \delta_j$ are pairwise distinct modulo the center of $G$.  Choosing $g$ and setting $g_j := \delta_j g \delta_j^{-1}$, it suffices to prove that the conjugates of $g$ by all $\gamma_i \delta_j$ generate a free subgroup.  Thus, replacing $r$ by $rd$, we may assume $d=1$.  We write $g$ for $g_1$.

Our main theorem is as follows.

\begin{thm}
\label{main}
Conjecture~\ref{indep} holds when $\uG$ has root system of type $A$, $B$, $E_6$, or $G_2$.
\end{thm}

The condition on $g$ that the elements $\gamma_i g\gamma_i^{-1}$
satisfy no relation is a countable union of conditions, each of which is closed in the Zariski-topology.  The set of $\R$-points of any proper closed subvariety of $\uG$ is a closed subset of $G$ in the real topology which has no interior points, so by the Baire category theorem, it suffices to prove that there is no non-trivial word (without coefficients) in the $\gamma_i g\gamma_i^{-1}$ which holds identically.

So Theorem~\ref{main} follows from the following variant.

\begin{thm}
\label{complex}
Let $\uG$ be a simple algebraic group over $\C$ of type $A$, $B$, $E_6$, or $G_2$.  Let $\gamma_1,\ldots,\gamma_r$ denote elements of $\uG(\C)$ such that every element of the form $\gamma_j^{-1}\gamma_i$, $i\neq j$, is semisimple but not central.  Then there exists $g\in \uG(\C)$ such that $\{\gamma_i g \gamma_i^{-1}\mid 1\le i\le r\}$ generates a free group of rank $r$.
\end{thm}

The rest of the paper is devoted to a proof of this theorem. 
Note that the claimed property of $\uG$ is invariant under isogenies, so we may choose $\uG$ freely within its isogeny class.  We use the simply connected form except for type $B_r$, where we use $\SO_{2r+1}$. 
We make use of a faithful representation of complex algebraic groups $\rho\colon \uG(\C)\to \GL(V)$.  For types $A_r$ and $B_r$, we use the natural representation.  For types $E_6$ and $G_2$ we use irreducible representations
of degree $27$ and $7$ respectively.

For any particular $\uG$, using the Baire category theorem as before, it suffices to prove
that there exists an extension field $L/\C$ and an element
$g\in \uG(L)$ such that $\{\gamma_i g \gamma_i^{-1}\mid i\in [1,r]\}$ generate a free subgroup of $\uG(L)$.  
We  set $L = \C((t))$ and define the valuation $v$ on $L$ so $t^{-v(a)}a\in \C[[t]]^\times$ for all $a\in L^\times$.
An $n+1$-tuple $(a_0,\ldots,a_n)\in L^{n+1}$ representing a point $[a_0:\cdots:a_n]\in L\P^n$ is said to be \emph{normalized} if $a_i\in \C[[t]]$ for all $i$ and $\sup_i v(a_i) =0$.
We define the map $\pi\colon L\P^n\to \C\P^n$ sending
$$[a_0:\cdots:a_n]\mapsto [a_0(0):\cdots:a_n(0)]$$
whenever $(a_0,\ldots,a_n)$ is normalized.  

For any non-zero column vector $v\in V= \C^{n+1}$, we define $B_v\subset L\P^n$ to be the inverse image by $\pi$ of the point in $\C\P^n$ determined by $v$.
For any non-zero row vector $v^*\in V^*$, we define $H_{v^*}$ to be the inverse 
image by $\pi$ of the hyperplane in $\C\P^n$ determined by $v^*$.
Let $e_0,\ldots,e_n$ denote the standard basis of $V$ and $e_0^*,\ldots,e_n^*$
its dual basis.

If $k_0<k_1\le \ldots\le k_{n-1} <k_n$ is an integer sequence, we define $\tau = t(k_0,\ldots,k_n)\in \GL_{n+1}(L)$ to be the diagonal matrix with entries $t^{k_i}$.
Thus, $\tau$ maps the complement of $H_{e_0*}$ to $B_{e_0}$ while $\tau^{-1}$ maps the complement of $H_{e_n^*}$ to $B_{e_n}$.

More generally, for $\eta\in \GL_{n+1}(\C)$, we have
\begin{equation}
\label{plus}
\eta \tau\eta^{-1}\colon H_{\eta(e_0^*)}^c \mapsto B_{\eta (e_0)}
\end{equation}
and 
\begin{equation}
\label{minus}
\eta \tau^{-1}\eta^{-1} \colon H_{\eta(e_n^*)}^c \mapsto B_{\eta (e_n)}.
\end{equation}
The following proposition gives the independence statement which is the key step in proving Theorem~\ref{complex}.

\begin{prop}
\label{pingpong}
Let  $\tau$ be defined as above, and let $\eta_1,\ldots,\eta_r$ be elements of $\GL_{n+1}(\C)$.
Suppose that for $i\neq j$ integers in $[1,r]$, we have
\begin{equation}
\label{nonincidence}
\eta_i(e_0), \eta_i(e_n)\not\in  \ker \eta_j(e^*_0) \cup \ker \eta_j(e^*_n).
\end{equation}
Then the set
$$\{g_i:=\eta_i \tau \eta_i^{-1}\mid i=1,2,\ldots,r\}$$
freely generates a free subgroup of $\GL_{n+1}(L)$.
\end{prop}

\begin{proof}
Applying $\pi^{-1}$, we can express \eqref{nonincidence} in the equivalent form 
\begin{equation}
\label{clean nonincidence}
B_{\eta_i(e_0)} \cup B_{\eta_i(e_n)} \subset H^c_{\eta_j(e^*_0)} \cap
H^c_{\eta_j(e^*_n)}
\end{equation}
for $i\neq j$.  We also have $e^*_i(e_i)\neq 0$, so
\begin{equation}
\label{repeated letter}
B_{\eta_i(e_0)} \subset H^c_{\eta_i(e^*_0)},\;
B_{\eta_i(e_ n)} \subset H^c_{\eta_i(e^*_n)}.
\end{equation}

There exists $z\in \C\P^n$ which does not lie in the hyperplane determined by any $\eta_i(e_0^*)$ or $\eta_i(e_n^*)$.
Let $\tilde z$ denote any element of $\pi^{-1}(z)$, so for all $i$,
\begin{equation}
\label{z-tilde}
\tilde z\in H^c_{\eta_i(e^*_0)} \cap H^c_{\eta_i(e^*_n)}.
\end{equation}
Let $x_{i_1}^{\epsilon_1}\cdots x_{i_l}^{\epsilon_l}$, $\epsilon_j \in \{\pm1\}$, $i_j\in [1,r]$ denote any non-trivial reduced word in $x_1,\ldots,x_r$.  
Starting at $j=l$, we prove statements \eqref{in ball} and \eqref{out of hyperplane} below by descending induction:
\begin{equation}
\label{in ball}
g_{i_j}^{\epsilon_j}\cdots g_{i_l}^{\epsilon_l}\tilde z \in
\begin{cases}
B_{\eta_{i_j}(e_0)}&\text{ if }\epsilon_j=1 \\
B_{\eta_{i_j}(e_n)}&\text{ if }\epsilon_j=-1.\\
\end{cases}
\end{equation}
\begin{equation}
\label{out of hyperplane}
g_{i_j}^{\epsilon_j}\cdots g_{i_l}^{\epsilon_l}\tilde z \in
\begin{cases}
H_{\eta_{i_{j+1}}(e^*_0)}^c&\text{ if }\epsilon_{j+1}=1 \\
H_{\eta_{i_{j+1}}(e^*_n)}^c&\text{ if }\epsilon_{j+1}=-1.\\
\end{cases}
\end{equation}
Note that if $j=l$ or $i_j \neq i_{j+1}$, then by \eqref{z-tilde} or \eqref{clean nonincidence} respectively, condition \eqref{in ball} implies condition \eqref{out of hyperplane}.  If $j<l$ and $i_j=i_{j+1}$, then $\epsilon_j = \epsilon_{j+1}$, so by \eqref{repeated letter}, 
condition \eqref{in ball} implies condition \eqref{out of hyperplane}.  Either way, for the induction step, we need only verify condition \eqref{in ball},  and this, in turn, follows from  (\ref{plus}) and (\ref{minus}).

Applying the induction hypothesis in the case $j=1$, we see that $\tilde z$ is not a fixed point of $g_{i_1}^{\epsilon_1}\cdots g_{i_l}^{\epsilon_l}$, so $g_{i_1}^{\epsilon_1}\cdots g_{i_l}^{\epsilon_l}\neq 1$.
\end{proof}

Let $\rho\colon \uG\to\GL_{n+1}$ denote a homomorphism of algebraic groups over $\C$ giving
rise to an irreducible representation of $\uG(\C)$.  Let
$\tau\in \uG(L)$ satisfy $\rho(\tau) = t(k_0,\ldots,k_n)$.  Suppose 
that there exists $h\in \uG(\C)$ such that setting $\eta_i = \rho(h^{-1}\gamma_i h)$, condition \eqref{nonincidence} holds for the $\eta_i$.
Then by Proposition~\ref{pingpong}, the set $\{\eta_i \rho(\tau)\eta_i^{-1}\}$ satisfies no relations in $\GL_{n+1}(L)$.
So setting $g = h \tau h^{-1}$, the set $\{\gamma_i g \gamma_i^{-1}\}$ satisfies no relations in $\uG(L)$.  This implies Theorem~\ref{complex} for $\uG$.

We begin with the following lemma:

\begin{lem}
\label{pairing}
Let $\uG$ be a semisimple algebraic group over $\C$, and let $\rho\colon \uG(\C)\to \GL(V)$ denote an irreducible representation such that $\End\,V = V^*\otimes V$ is a direct sum of pairwise distinct irreducibles.
Let $T$ be an element of $\End\,V$.
Suppose that there exists an irreducible factor $W$ of $\End\,V$ such that $T$ has a non-zero component in $W$ and $v^*\otimes v$ has non-zero component in $W^*$.  Then
$$\{h\in \uG(\C)\mid \rho(h)(v^*)(T\rho(h) v)\neq 0\}$$
is a dense Zariski-open subset of $\uG(\C)$.
\end{lem}

\begin{proof}
The openness is clear, so it suffices to prove that there exists $h$ with non-vanishing
$\rho(h)(v^*)(T\rho(h) v)$.
This expression can be rewritten 
$$\langle\rho(h)(v^*\otimes v),T\rangle,$$
where $\langle\;,\,\rangle$ denotes the usual pairing $\End\,V\times \End\,V\to \C$.
If it is identically zero, then the pairing of $T$ with every element of
\begin{equation}
\label{spanning set}
\Span\{\rho(h)(v^*\otimes v)\mid h\in \uG(\C)\}
\end{equation}
is zero.  By the multiplicity one property of $\End\,V$, $\eqref{spanning set}$ is the sum of all factors in which the projection of $v^*\otimes v$ is non-zero, which by assumption includes the factor $W^*$.
As $\langle\;,\,\rangle$ is a perfect pairing, it induces a perfect pairing between $W$ and $W^*$, so the pairing of $T$ with a suitable element of $W^*$ is non-zero, and the lemma holds.
\end{proof}

\begin{lem}
\label{Minuscule}
Let $\uG$ be a simple algebraic group over $\C$, and let $V$ be a non-trivial minuscule representation of $\uG(\C)$.  Let $\uT\subset \uB$ be a maximal torus and a Borel subgroup of $\uG$.  Let $\chi$ denote the maximal weight of $\uT$ acting on $V$, and let $v$ and $v^*$ denote non-zero vectors in $V$ and $V^*$ respectively, in the weight spaces $\chi$ and $-\chi$ respectively.  Then 
\eqref{spanning set} is $\End \,V$.
\end{lem}

\begin{proof}
As the Weyl group $N(\uT)/\uN$ acts transitively on the weights of $V$ with respect to $\uT$, the span
\begin{equation}
\label{N-span}
\Span \{\rho(h)(v^*\otimes v)\mid h\in \uN(\C)\},
\end{equation}
and therefore the span \eqref{spanning set}, contains a non-zero element of $V^*\otimes v'$ for every non-zero vector $v'$ belonging to any weight space of $V$.  Therefore, the dimension of \eqref{N-span} is at least $\dim V$.  On the other hand, every vector in this space lies in the trivial $\uT$-weight space of $\End\,V$, which has dimension $\dim V$.  Therefore,
the span of \eqref{spanning set} contains the trivial weight space of $\End\,V$.

On the other hand, every irreducible factor $W$ of $\End\,V$ has trivial central character as $\uG$-representation and therefore has highest weight in the root lattice; it follows that $0$ is a weight of each such $W$.  Therefore, the image of \eqref{spanning set} under projection to $W$ contains a non-zero vector, and since it is $\uG$-stable, it must span $W$.  As the image of the span of \eqref{spanning set} is a subrepresentation of $\End\,V$ which maps onto every irreducible quotient representation of $\End\,V$,
it must be all of $\End\,V$.

\end{proof}

We recall that an irreducible representation of a simple algebraic group is \emph{quasi-minuscule} if the non-zero weights form a single orbit under the Weyl group.  For instance, the natural representation of $\SO_{2k+1}$ is quasi-minuscule for $k\ge 2$.  Among other examples, the adjoint representation is quasi-minuscule for simply laced root systems, and $G_2$ has a quasi-minuscule representation of degree $7$.
\begin{prop}
\label{BG}

For $\uG = \SO_{2k+1}$, $k\ge 2$, and $V$ the natural representation of $\uG(\C)$, if $v$ and $v^*$ denote, respectively, non-zero vectors in the highest weight space of $V$ and the lowest weight space of $V^*$, then \eqref{spanning set} contains $\sl(V)$.
The same is true for the $7$-dimensional representation of $G_2$.
\end{prop}

\begin{proof}
For $\SO_{2k+1}$, preceding as before, we see that the dimension of \eqref{N-span}
%
is $2k$.  Identifying $V^*$ with $V$, $\End\,V\cong V^{\otimes 2}$ decomposes into the symmetric square and the exterior square, and both the symmetric and the anti-symmetric parts of $\sl(V)$ are irreducible.  Each has a $k$-dimensional trivial weight space, while the remaining (trivial) factor of $\End\,V$ has a $1$-dimensional trivial weight space.  Therefore, the image of \eqref{N-span} in each of these irreducible factors is non-zero, and we finish as before.

When $V$ is the $7$-dimensional irreducible representation of $G_2$ the dimension of \eqref{N-span} is $6$, and its intersection with $\Sym^2 V$ and $\wedge^2 V$ are both $3$-dimensional.
By \cite[Table~A.121]{FKS}, the trace-$0$ subrepresentation  $\sl(V)\subset \End\,V$
decomposes into irreducible factors of degrees $7$, $14$, and $27$.  For dimension reasons, it follows that
$\wedge^2 V\subset \sl(V)$ decomposes into the $7$ and the $14$, while
$\Sym^2 V \cap \sl(V)$ is the 27.  Each has $3$-dimensional trivial weight space, which therefore coincides with its intersection with \eqref{N-span}.
Since the weight and root lattices are the same for $G_2$, every non-trivial representation has 
positive dimensional trivial weight space, and it follows that
the image of \eqref{N-span} in each of the three irreducible factors has positive dimension.  We finish as before.
\end{proof}

\begin{prop}
\label{Diagonal}
Let $\rho\colon \uG(\C) \to \GL(V)$ be an irreducible orthogonal representation of degree $2k+1$ such that $\Sym^2 V\cap \sl(V)$
is irreducible.
Let $v$ be a non-zero vector in $V$.  Then for all $\gamma\in \uG(\C)$, there exists $h\in \uG(\C)$ such that 
$ \rho(h) (v)\cdot  \rho(\gamma h) (v) \neq 0$.
\end{prop}

\begin{proof}
We identify $V$ with $V^*$ and therefore $V\otimes V$ with $\End\,V$.
As an endomorphism of $V$, $v\otimes v$ has rank exactly $1$, so it is not a scalar multiple of the identity.  Thus, its projection onto $\Sym^2 V\cap \sl(V)$ is non-zero, so the span of its $\uG(\C)$-orbit projects onto $\Sym^2 V\cap \sl(V)$. If $ \rho(h)(v)\cdot  \rho(\gamma h)(v)= 0$ for all $h\in \uG(\C)$, then $\rho(\gamma)\in \rho(\uG(\C))\subset \End\,V$ must pair to $0$ with each $\rho(h)v\otimes \rho(h)v$ and therefore with every element in $\Sym^2 V$.  It follows that it must lie in the sum of the scalar matrices and $\wedge^2 V$.  

We claim that a diagonalizable orthogonal matrix $O\in \SO_{2k+1}(\C)$, $O\neq I$, cannot be the sum $cI+M$ of a scalar matrix and a skew-symmetric matrix.
As the $\lambda$-eigenspace of $O$ is orthogonal to the $\mu$-eigenspace unless $\lambda\mu=1$, there exists a basis $f_0,\ldots,f_{2n}$ of $V$ consisting of eigenvectors and such that
$\Span(f_0),\,\Span(f_1,f_2),\,\Span(f_3,f_4),\,\ldots,\,\Span(f_{2k-1},f_{2k})$ are mutually orthogonal.  As each of these subspaces has an orthonormal basis, after conjugating $O$ by a suitable orthogonal matrix
we may assume that it is in block diagonal form, with blocks of size $1,2,2,\ldots,2$, each of which is orthogonal and such that $O$ is still of the form $cI + M$, with $M$ skew-symmetric.  As $\det(O)=1$, $O f_0 = f_0$, while $Mf_0 = 0$, so
$c=1$.  However, the identity matrix itself is the only orthogonal $2\times 2$ matrix which is the identity plus a skew-symmetric matrix, thus proving our claim.  The claim implies the proposition.
\end{proof}

We remark that this applies to the natural representation $V$ when $\uG$ is an odd orthogonal group and to the $7$-dimensional representation when $\uG$ is of type $G_2$.

\begin{prop}
\label{e6}
Let $\uG$ denote the simply connected simple algebraic group of type $E_6$ over $\C$ and let $\rho\colon \uG(\C)\to \GL(V)$ denote a minuscule irreducible representation.
Let $\rho_*\colon \g\to \End\,V$ denote the corresponding Lie algebra representation.  If $1\neq h\in \uG(\C)$, then $\rho(h)$ cannot be expressed in $\End\,V$ as the sum of a scalar matrix and a non-zero element of $\rho_*(\g)$.
\end{prop}

\begin{proof}
In \cite{Larsen}, it is proved that for exceptional complex matrix groups, $\rho(h)$ cannot lie in $\rho_*(\g)$, and the argument here is nearly identical to the $E_6$ case.  Suppose that $\rho(h) = \rho_*(x) + cI$.
We fix a maximal torus $\uT$ containing $h$ and a basis of $V$ consisting of weight vectors with respect to $\uT$.  Let $\frakt$ denote the Lie algebra of $\uT(\C)$, so $\frakt$ contains $x$.
If $\chi_i,\chi_j,\chi_k,\chi_l$ are four characters of $\uT$  appearing as weights of $V$ and $t^*_i,t^*_j,t^*_k,t^*_l$  the corresponding elements of $\frakt^*$,
then $\chi_i \chi_j = \chi_k\chi_l$ implies
$$\chi_i(g)+\chi_j(g) = t^*_i(x) + c + t^*_j(x)+c = t^*_k(x) + c + t^*_l(x)+c = \chi_k(g)+\chi_l(g).$$
Thus, \cite[equation (1)]{Larsen} holds, and this is all that is needed to prove the version \cite[Proposition~6]{Larsen} in which $\rho_*(\g)$ is replaced by $\rho_*(\g)+\C\,I$.
\end{proof}

\begin{prop}
Let $\uG$ be a simple algebraic group over $\C$, $\rho\colon \uG\to \GL(V)$ an irreducible representation, and  $\gamma_1,\ldots,\gamma_r\in \uG(\C)$  pairwise distinct modulo the center of $\uG(\C)$.  In each of the following situations, there exists $h\in \uG(\C)$ such that the sequence 
$$\eta_1=\rho(h^{-1} \gamma_1 h),\,\ldots,\,\eta_r = \rho(h^{-1} \gamma_r h)$$
satisfies \eqref{nonincidence}:
\begin{enumerate}
\item[1.]$\uG = \SL_n$ and $V$ is the natural representation.
\item[2.]$\uG = \SO_{2k+1}$ and $V$ is the natural representation.
\item[3.]$\uG = G_2$ and $\dim V = 7$.
\item[4.]$\uG$ of type $E_6$ and $\dim V=27$.
\end{enumerate}

\end{prop}

\begin{proof}
Let $\uT$ be a maximal torus of $\uG$.  Conjugating $\rho$ if necessary, we may assume $\rho(\uT)$ is contained in the set of diagonal matrices of $\GL_{n+1}(\C) = \GL(V)$.  Let $c\colon \G_m\to \uT$ be a cocharacter with the property that its composition with the weights of $\rho$ with respect to $\uT$ are pairwise distinct.  Then $\tau=\rho(c(t))$ is a matrix of the form $t(k_0,\ldots,k_n)$ with pairwise distinct $k_i$, and conjugating $\rho$ if necessary, we may assume the sequence $k_0,\ldots,k_n$ is strictly increasing.

Let $v$ be a non-zero vector in the highest weight space of $V$ as $\uG(\C)$-representation.
Let $v^*$ and $w^*$ denote, respectively,
a non-zero vector in the lowest (resp. the highest) weight space of $V^*$.  The condition 
$$\eta_i(v)\not\in \ker \eta_j (v^*)$$
from \eqref{nonincidence} can therefore be written
$$\rho(\gamma_j^{-1}\gamma_i) \rho(h)B_v \subset \rho(h) H^c_{v^*}$$
or
$$\langle \rho(h)(v^*\otimes v),\rho(\gamma_j^{-1}\gamma_i)\rangle =  (\rho(g)v^*)(\rho(\gamma_j^{-1}\gamma_i)(\rho(g) v)) \neq 0.$$
for some $g$.  Likewise, the condition
$$\eta_i(v)\not\in \ker \eta_j (w^*)$$
can be written
$$\langle \rho(h)(w^*\otimes v),\rho(\gamma_j^{-1}\gamma_i)\rangle \neq 0.$$

In case 1, $\End\,V = \sl(V)\oplus \C$ is the decomposition into irreducible factors.  In case 2, it is $\End\,V = \wedge^2 V \oplus (\Sym^2 V\cap \sl(V)) \oplus \C$.
In case 3, we have seen that the irreducible factors have degrees $1$, $7$, $14$, and $27$.  In case 4, $\End\,V$ has one trivial factor, one $78$-dimensional factor which is the adjoint representation, and one $650$-dimensional factor \cite[Table A.117]{FKS}.
In none of these cases is there a repeating factor, so Lemma~\ref{pairing} applies.

In case 1, by hypothesis, $\rho(\gamma_j^{-1}\gamma_i) = \gamma_j^{-1}\gamma_i$ is not scalar, and $v^*\otimes v$ and $w^*\otimes v$ have rank $1$, so they are not scalar either.  Thus, all three of these vector have non-trivial projection to the irreducible factor $\sl(V)$, and we are done.
In cases 2 and 3, Proposition~\ref{BG} covers $v^*\otimes v$ and Proposition~\ref{Diagonal} covers $w^*\otimes v$.
In case 4, the projections of $\rho(\gamma_j^{-1}\gamma_i)$ and $v^*\otimes v$ onto the $650$-dimensional factor are non-trivial by Proposition~\ref{e6} and Lemma~\ref{Minuscule} respectively.
The vector $w^*\otimes v$ lies in the highest weight space of the $650$-dimensional space, so we are done.

\end{proof}

\end{document}